\documentstyle{amsart}
\title[Nonlinear Algebraic Systems...]
{Nonlinear Algebraic Systems\\ with discontinuous terms}

\author{Giovanni Molica Bisci and Du\v{s}an Repov\v{s}}

\address[G. Molica Bisci]{Dipartimento MECMAT, University of Reggio Calabria,
 Via Graziella, Feo di Vito, 89124 Reggio Calabria, Italy.} \email{gmolica@@unirc.it}
\address[D. Repov\v{s}]{Faculty of Education, and Faculty of Mathematics and Physics, University of Ljubljana, Kardeljeva pl. 16, Ljubljana, Slovenia 1000.} \email{dusan.repovs@@guest.arnes.si}

\thanks{}



\newtheorem{theorem}{Theorem}[section]
\newtheorem{proposition}{Proposition}[section]
\newtheorem{corollary}{Corollary}[section]
\newtheorem{remark}{Remark}[section]

\def\essinf{\mathop {\rm ess\,inf}}
\def\esssup{\mathop {\rm ess\,sup}}

\def\erre{{\rm I\!R}}
\def\essinf{\mathop {\rm ess\,inf}}
\def\esssup{\mathop {\rm ess\,sup}}
\def\erre{{\rm I\!R}}
\def\R{{\rm I\!R}}

\def\essinf{\mathop {\rm ess\,inf}}

\def\phi{\varphi}

\def\trace{\mathop{\rm Tr}(A)}
\def\T{\mathop{\rm Trid}}
\begin{document}

\begin{abstract}
Using a multiple critical points theorem for locally Lipschitz continuous functionals, we establish the existence of at least three distinct solutions for a parametric discrete differential inclusion problem involving a real symmetric and positive definite matrix. Applications to tridiagonal, fourth-order and partial difference inclusions are presented.
\end{abstract}

\maketitle

\section{Introduction}
Let $T>1$ be a positive integer and let $g_k:\erre\rightarrow\erre$ be locally essentially bounded functions,
for every $k\in {\mathbb{Z}}[1,T]:=\{1,2,...,T\}$. In this paper we are interested in the existence of multiple solutions
for the following discrete problem
\begin{equation} \tag{$S_{A,\lambda}^{g}$} \label{W}
\sum_{l=1}^{T}a_{kl}u_l\in \lambda[g^-_k(u_k),g^+_k(u_k)],\quad\quad \forall\,k\in {\mathbb{Z}}[1,T],
\end{equation}
where $\lambda$ is a positive parameter, $A:=(a_{ij})_{T\times T}$ is a real symmetric positive definite matrix and
$$
g^{-}_k(t):=\lim_{\delta\to 0^{+}}\essinf_{|\xi-t|<\delta}g_k(\xi),\,\,\,\,\,\,\,\,\,
g^{+}_k(t):=\lim_{\delta\to 0^{+}}\esssup_{|\xi-t|<\delta}g_k(\xi),
$$
for every $k\in {\mathbb{Z}}[1,T]$. It is clear that if the functions $g_k$ are continuous (instead of locally essentially bounded) problem \eqref{W} becomes a more familiar nonlinear algebraic system
$$
Au=\lambda g(u),
$$
 in which $u=(u_1,...,u_T)^{t}\in\erre^T$ and $g(u):=(g_1(u_1),...,g_T(u_T))^t$.\par

A considerable number of problems, which are strictly connected both with boundary value differential problems and numerical simulations of some mathematical models arising from many research areas (e.g. biological, physical and computer science) can be formulated as special cases of nonlinear algebraic systems.\par
However, to the
best of our knowledge, for discrete difference inclusions there are only few papers involving
the second-order difference operator. For instance, in \cite{AOL}, the existence of at least
one solution was obtained via the set-valued mapping theory, while in \cite{Z45}, existence results
for suitable second-order discrete discontinuous equations have been investigated by variational methods.\par
The aim of this paper is to establish a three solutions result (Theorem \ref{Main} below) for problem \eqref{W}. Successively, in Corollary \ref{Main2}, the existence of at least two nontrivial solutions is obtained for sufficiently large values of the real parameter $\lambda$.\par
 Our proof is based on a three critical points theorem
for nondifferentiable functionals obtained by Bonanno and Marano (see Theorem \ref{abstract} in Section 2).
We also note that, very recently, Theorem \ref{abstract} has been used in \cite{CL}, in order to obtain an analogous version of Theorem \ref{Main} for Dirichlet inclusions involving the discrete $p$-Laplacian operator.\par
Due to the generality of problem \eqref{W} remarkable applications are achieved.
For instance, in Section 4 we present some existence theorems
concerning discrete inclusions involving certain tridiagonal matrices, fourth-order discrete inclusions and partial difference inclusions.\par
A special case of Corollary \ref{Main2} reads as follows.

\begin{theorem}\label{Introduction}
Let $h:\erre\to\erre$ be essentially locally bounded positive function and consider the usual forward difference operator $\Delta
u_{k-1}:=u_{k}-u_{k-1}$.
Assume that
\begin{itemize}
\item[$(\rm h_2)$] $\displaystyle\limsup_{t\rightarrow 0^+}\frac{h(t)}{t}=0$$;$
\item[$(\rm h_3)$] $\displaystyle\limsup_{|t|\rightarrow \infty}\frac{ h(t)}{t}<\displaystyle\frac{\lambda_1}{T},$ where $\lambda_1$ is the first eigenvalue of the discrete problem
$$
\left\{
\begin{array}{l}
\Delta^4u_{k-2}=\lambda u_k,\quad\forall\, k  \in
{\mathbb{Z}}[1,T] \\
{u_{-2}=u_{-1}=u_0=0}\\
u_{T+1}=u_{T+2}=u_{T+3}=0.\\
\end{array}
\right.
$$
\end{itemize}
\noindent Then for every
$$
\lambda>\frac{2}{T}\inf_{\delta>0}\frac{\delta^2}{\displaystyle\int_{0}^{\delta}h(t)dt},
$$
the following forth-order discrete inclusion
\begin{equation} \tag{$D_{\lambda}^{h}$} \label{4order}
\left\{
\begin{array}{l}
\Delta^4u_{k-2}\in \lambda [h^-(u_k),h^+(u_k)],\quad \forall\,k\in {\mathbb{Z}}[1,T]\\
{u_{-2}=u_{-1}=u_0=0}\\
u_{T+1}=u_{T+2}=u_{T+3}=0.\\
\end{array}
\right.
\end{equation}
 admits at least two nontrivial solutions.
\end{theorem}

Finally, we emphasize that our results are also new in the continuous setting. In this case, the existence and multiplicity of
solutions was investigated in a large number of other papers under various assumptions (see
for instance \cite{Z,ZB,ZC,ZF,Z45} and references therein). See also recent papers \cite{CMolica, KMR, KMRT, MRT} for related topics.\par
The plan of the paper is as follows. In the next section we introduce our abstract framework. Successively, in section 3, we show our multiplicity results.\par In conclusion, concrete examples of applications of our abstract results are presented.
\section{Abstract framework}
\indent Let $(X,\Vert\cdot\Vert)$ be a real Banach space.
We denote by $X^*$ the dual space of $X$, while
$\langle\cdot,\cdot\rangle$ stands for the duality pairing between
$X^*$ and $X$. A function $J:X\to\erre$ is called locally
Lipschitz continuous if to every $x\in X$ there corresponds a
neighborhood $V_x$ of $x$ and a constant $L_x\geq 0$ such that
$$\vert J(z)-J(w)\vert\leq L_x\Vert z-w\Vert,\quad\forall\, z,w\in V_x\, .$$
\indent If $x,z\in X$, we write $J^{0}(x;z)$ for the generalized
directional derivative of $J$ at the point $x$ along the
direction $z$, i.e.,
$$J^{0}(x;z):=\limsup_{w\to x,\, t\to 0^+}{{J(w+tz)-J(w)}\over t}\, .$$
\noindent The generalized gradient of the function $J$ in $x$,
denoted by $\partial J(x)$, is the set
$$\partial J(x):=\left\{ x^*\in X^*:\, \langle x^*,z\rangle\leq J^{0}(x;z),
\;\forall\, z\in X\right\}.$$
 \indent The basic properties of generalized directional derivative and generalized gradient were studied in \cite{CHANG,CLARKE2}.
We recall that if $J$ is continuously G\^{a}teaux differentiable at $u$, then $J$ is locally
Lipschitz at $u$ and $\partial J(u) = \{J '(u)\}$, where $J'(u)$ stands for the first derivative of $J$ at $u$.\par
Further, a point $u$ is called a (generalized)
critical point of the locally Lipschitz continuous function $J$ if $0_{X^*}\in \partial J(u)$,
 i.e.
$$J^0(x;z)\geq 0,$$
 for every $z\in X$. Clearly, if $J$ is a continuously G\^{a}teaux differentiable at $u$, then $u$
becomes a (classical) critical point of $J$, that is $J '(u) = 0_{X^*}$.\par
For an exhaustive overview on the non-smooth calculus we mention the excellent
monograph \cite{MR}. Further, we cite a very recent book \cite{KRV} as a general reference on this subject.\par
The main tool will be the following abstract critical point theorem for locally Lipschitz continuous functions on finite dimensional Banach spaces that can be derived from \cite[Theorem 3.6]{BONMAR}.

\begin{theorem}\label{abstract}
{Let $X$ be a finite dimensional real Banach space, $\Phi, \Psi :X\to\R$ be two locally Lipschitz continuous functionals such that $\Phi(0_X)=\Psi(0_X)=0$ and
assume that the following conditions are satisfied$:$
\begin{itemize}
\item[$({\rm a}_1)$] there exist $r>0$ and $\bar{u}\in X$, with $\Phi(\bar{u})>r$, such that$:$
$$
\displaystyle\frac{\displaystyle\sup_{u \in \Phi^{-1}(]-\infty,r])}\Psi(u)}{r}< \frac{\Psi(\bar{u})}{\Phi(\bar u)};
$$
\item[$({\rm a}_2)$] for each $\lambda\in \Lambda_{r}:=\left]\displaystyle\frac{\Phi(\bar u)}{\Psi(\bar{u})},\displaystyle\frac{r}{\displaystyle\sup_{u \in \Phi^{-1}(]-\infty,r])}\Psi(u)}\right[$ the functional $J_{\lambda}:=\Phi-\lambda\Psi$ is coercive, that is
    $$
    \lim_{\|u\|\rightarrow \infty}(\Phi(u)-\lambda\Psi(u))=+\infty.
    $$
\end{itemize}
Then for each $\lambda \in \Lambda_{r}$, the functional $J_{\lambda}$ has at least three distinct critical points in $X$.
}
\end{theorem}

We also cite the quoted papers of Ricceri \cite{R1,R2,R3} for some related topics in the smooth context.\par
\smallskip

\indent Here, as the ambient space $X$, we consider the $T$-dimensional Banach space $\erre^T$
endowed by the norm
$$
\|u\|:=\Big(\sum_{k=1}^{T}u_k^2\Big)^{1/2},
$$
induced by the standard Euclidean inner product $\langle u,v\rangle_X:=\displaystyle\sum_{k=1}^{T}u_kv_k$.\par
\indent Set ${\mathfrak{X}}_{T}$ to be the class of all symmetric and positive definite matrices of order $T$. Further, we denote by $\lambda_1,...,\lambda_T$ the eigenvalues of $A$ ordered as follows $0<\lambda_1\leq...\leq \lambda_T$.\par
It is well-known that if $A\in {\mathfrak{X}}_{T}$, for every $u\in X$, then one has

\begin{equation}\label{immersione2}
\lambda_1\|u\|^2\leq u^{t}Au\leq \lambda_T\|u\|^2,
\end{equation}
\noindent and
\begin{equation}\label{immersione}
\|u\|_\infty\leq \frac{1}{\sqrt{\lambda_1}}(u^{t}Au)^{1/2},
\end{equation}
\noindent where $ \|u\|_{\infty}:=\displaystyle\max_{k\in
{\mathbb{Z}}[1,T]}|u_k| $.

\indent For every $u\in X$, put
\begin{equation} \label{functions}
\Phi(u):=\frac{u^{t}Au}{2}, \quad  \quad \Psi(u):=\sum_{k=1}^{T}G_k(u_k), \quad \textrm{and} \quad J_\lambda(u):=\Phi(u)-\lambda\Psi(u),
\end{equation}
\noindent where
$G_k(t):=\displaystyle\int_{0}^{t}g_k(\xi)d \xi$, for every $(k,t)
\in {\mathbb{Z}}[1,T]\times \erre$.\par

It is easy to verify that $\Phi$ is continuously G\^{a}teaux differentiable, while $\Psi$ is locally
Lipschitz continuous.

\begin{proposition}\label{equivalen}
Fix $\lambda>0$ and assume that $u\in X$ is a critical point of the functional $J_\lambda:=\Phi-\lambda\Psi$. Then $u$ is a solution of problem \eqref{W}.
\end{proposition}

\begin{pf}
If $u$ is a critical point of $J_\lambda$, bearing in mind of \cite[Propositions 2.3.1 and
2.3.3]{CLARKE2}, it follows that
\begin{equation} \label{derivata}
\Phi'(u)(z)\leq \Psi^{0}(u;z)\leq \lambda\left(\sum_{k=1}^{T}G^{0}_k(u_k;z_k)\right),
\end{equation}
for every $z\in X$. Moreover,
\begin{equation} \label{rappresentazione}
\Phi'(u)(z)=\frac{\langle\nabla (u^{t}Au),z\rangle_{X}}{2},
\end{equation}
for every $z\in X$.
For every $\xi\in\erre$ and $k\in {\mathbb{Z}}[1,T]$, by putting in \eqref{derivata} the vector $z=\xi e_k$,
where $e_k$ are the canonical unit vectors of $X$, and taking in mind \eqref{rappresentazione}, we get
$$
\langle \sum_{l=1}^{T}a_{kl}u_l,\xi\rangle_{\erre}=\Phi'(u)(z)\leq \lambda G^{0}_k(u_k;\xi),
$$
namely
$$
\sum_{l=1}^{T}a_{kl}u_l\in \lambda \partial G_k(u_k).
$$
Finally, since it is well-known that
$$
\partial G_k(u_k)=[g^-_k(u_k),g^+_k(u_k)],
$$
for every $k\in {\mathbb{Z}}[1,T]$ (see for instance \cite[Example 2.2.5]{CLARKE2}) it follows that
$$
\sum_{l=1}^{T}a_{kl}u_l\in \lambda[g^-_k(u_k),g^+_k(u_k)],\quad\quad \forall\,k\in {\mathbb{Z}}[1,T].
$$
Therefore our assertion is proved.
\end{pf}
\section{Main Results}

In this section we present our multiplicity results for problem \eqref{W} which can be deduced by using Theorem \ref{abstract}.

\begin{theorem}\label{Main}
Let $g_k:\erre\rightarrow\erre$ be a locally essentially bounded function
for every $k\in {\mathbb{Z}}[1,T]$.
Assume that
\begin{itemize}
\item[$(\rm g_1)$] There exist positive constants $\gamma$ and $\delta$, with $$\delta>\displaystyle\left(\frac{\lambda_1}{\trace +2\sum_{i<j}a_{ij}}\right)^{1/2}\gamma,$$ such that
$$
\frac{\displaystyle\sum_{k=1}^{T}\sup_{|\xi|\leq \gamma} G_k(\xi)}{\gamma^2}<\frac{\lambda_1}{\trace +2\sum_{i<j}a_{ij}}\frac{\displaystyle\sum_{k=1}^{T}G_k(\delta)}{\delta^2};
$$
\end{itemize}

\begin{itemize}
\item[$(\rm g_2)$] $\displaystyle\limsup_{|\xi|\rightarrow \infty}\frac{\displaystyle G_k(\xi)}{\xi^2}<\frac{\lambda_1}{2},\,\,\,\,\,\,\forall\, k\in {\mathbb{Z}}[1,T].$

\end{itemize}
\noindent Then, for every $\lambda$ belonging to
$$
\Lambda:=\left]\frac{\trace +2\sum_{i<j}a_{ij}}{2}\frac{\delta^2}{\displaystyle\sum_{k=1}^{T} G_k(\delta)},\frac{\lambda_1}{2}\frac{\gamma^2}{\displaystyle\sum_{k=1}^{T}\sup_{|\xi|\leq \gamma} G_k(\xi)}\right[,
$$
problem \eqref{W} admits at least three solutions.
\end{theorem}
\begin{pf}
Fix $\lambda\in \Lambda$ and let
$\Phi$, $\Psi$,  $J_\lambda$ as
indicated in \eqref{functions}.
Since the critical points of $J_\lambda$ are the solutions of problem \eqref{W}, our aim is to apply Theorem \ref{abstract} to function $J_\lambda$.
Therefore, put
\[r:=\frac{\lambda_1}{2}\gamma^2,\]
\noindent and let us denote
$$
A_{k}:=\left\{u\in X: |u_k| \le \gamma\right\},$$
for every $k\in {\mathbb{Z}}[1,T]$.\par
\noindent Thanks to (\ref{immersione}), it follows that
$$
\{u\in X:u^{t}Au\leq 2r\}\subseteq \left\{u\in X: \|u\|_{\infty} \le \gamma\right\}\subseteq A_{k},
$$
for every $k\in {\mathbb{Z}}[1,T]$. Hence
$$
\chi(r):={\displaystyle\sup_{u \in \Phi^{-1}(]-\infty,r])}\Psi(u) \over
r}\leq\frac{\displaystyle\sup_{\|u\|_{\infty} \le \gamma}\sum_{k=1}^{T} G_k(u_k)}{r}\leq
$$
\begin{equation}\label{$3$}
\le\frac{\displaystyle\sum_{k=1}^{T}\sup_{u\in A_{k}} G_k(u_k)}{r}=\frac{2}{\lambda_1}\frac{\displaystyle\sum_{k=1}^{T}\sup_{|\xi|\leq \gamma} G_k(\xi)}{\gamma^2}.
\end{equation}
Now consider $\bar{u}\in X$ such that
$
\bar{u}_k:=\delta
$, for every $k\in {\mathbb{Z}}[1,T]$,
and observe that, since
$$
\delta>\displaystyle\left(\frac{\lambda_1}{\trace +2\sum_{i<j}a_{ij}}\right)^{1/2}\gamma,
$$
 one has
$
\Phi(\bar{u})>r.
$
Further,
\begin{equation}\label{$4$}
\frac{\Psi(\bar{u})}{\Phi(\bar{u})}=\frac{2}{\trace +2\sum_{i<j}a_{ij}}\frac{\displaystyle\sum_{k=1}^{T}G_k(\delta)}{\delta^2}.
\end{equation}
At this point, taking into account relations \eqref{$3$} and \eqref{$4$}, condition ($\rm g_2$) implies that
$$
\chi(r)\leq \frac{2}{\trace +2\sum_{i<j}a_{ij}}\frac{\displaystyle\sum_{k=1}^{T}G_k(\delta)}{\delta^2}.
$$
Therefore, hypothesis $({\rm a}_1)$ holds and $\Lambda\subseteq \Lambda_r$.\par
\noindent By ($\rm g_2$) there are constants $\epsilon \in \displaystyle\left]0,\lambda_1/2\right[$ and
$\sigma>0$ such that
\begin{equation}\label{z}
{\frac{\displaystyle\int_0^\xi g_k(t)dt}{\xi^{2}}}<\frac{\lambda_1}{2}-\epsilon,
\end{equation}
 for every
$|\xi|\geq \sigma$ and $k\in {\mathbb{Z}}[1,T]$. Let us put
\begin{equation}\label{x}
M_1:=\max_{(k,\xi)\in {\mathbb{Z}}[1,T]\times [-\sigma,\sigma]}\int_{0}^{\xi}g_k(t)dt.
\end{equation}
\noindent At this point note that, for every $\xi\in\erre$ and $k\in {\mathbb{Z}}[1,T]$, one has
$$
\int_0^\xi g_k(t)dt\leq M_1+M_2\xi^2,
$$
\noindent where
$$
M_2:=\frac{\lambda_1}{2}-\epsilon.$$

\noindent Moreover, the following inequality holds
$$
J_\lambda(u)\geq\frac{u^{t}Au}{2}-\sum_{k=1}^{T}\Big[M_1+M_2u_k^2\Big],\,\,\,\,\,\,\forall\, u\in X.
$$
\noindent Hence,
$$
J_\lambda(u)\geq \frac{u^{t}Au}{2}-M_2\|u\|^2-TM_1,\,\,\ \forall\, u\in X,
$$
 and, by relation (\ref{immersione2}), one has
\begin{eqnarray}\label{C}
J_\lambda(u)\geq \epsilon \|u\|^2-TM_1,\,\,\,\,\,
\forall\, u\in X,
\end{eqnarray}
\noindent which clearly shows that
\begin{eqnarray}\label{Coercivity}
\lim_{\|u\|\rightarrow\infty}J_\lambda(u)=+\infty.
\end{eqnarray}
So, the assumptions of Theorem \ref{abstract} are satisfied and our conclusions follow from
Proposition \ref{equivalen}.
\end{pf}
\begin{remark}
\rm{We point out that very recently, Theorem \ref{Main} has been exploited in \cite{CL}, in order to obtain an analogous version of Theorem \ref{Main} for Dirichlet inclusions involving the discrete $p$-Laplacian operator. A related variational approach has been also adopted, studying difference equations (in the continuous case) in papers \cite{BC1,BC2} as well as in \cite{CMolica2,IM} proving multiplicity results for nonlinear algebraic systems.}
\end{remark}

A direct consequence of Theorem \ref{Main} reads as follows.

\begin{corollary}\label{Main2}
Let $\alpha:{\mathbb{Z}}[1,T]\rightarrow\erre$ be a nonnegative $($not identically zero$)$ function and let $h:\erre\to\erre$ be a essentially locally bounded map.
Assume that
\begin{itemize}
\item[$(\rm h_1)$] there exists $\delta>0$ such that $h(t)>0$ for every $0<|t|<\delta$$;$
\item[$(\rm h_2)$] $\displaystyle\limsup_{t\rightarrow 0^+}\frac{h(t)}{t}=0$$;$
\item[$(\rm h_3)$] $\displaystyle\limsup_{|t|\rightarrow \infty}\frac{ h(t)}{t}<\displaystyle\frac{\lambda_1}{\displaystyle\sum_{k=1}^{T}\alpha_k}.$
\end{itemize}
\noindent Then for every
$$
\lambda>\frac{\displaystyle\trace/2 +\sum_{i<j}a_{ij}}{\displaystyle\sum_{k=1}^{T}\alpha_k}\frac{\delta^2}{\displaystyle\int_{0}^{\delta}h(t)dt},
$$
the following discrete problem
\begin{equation} \tag{$S_{A,\lambda}^{\alpha,h}$} \label{W2}
\sum_{l=1}^{T}a_{kl}u_l\in \lambda\alpha_k[h^-(u_k),h^+(u_k)],\quad\quad \forall\,k\in {\mathbb{Z}}[1,T],
\end{equation}
 admits at least two nontrivial solutions.
\end{corollary}
\begin{pf}\rm{
Our aim is to apply Theorem \ref{Main} to problem \eqref{W2}. Hence, let us
put $H(\xi):=\displaystyle\int_0^{\xi}h(t)dt$ for every $\xi\in\erre$. By $(\rm h_1)$ one has
$H(\delta)>0$ and
$$
\sup_{|\xi|\leq \gamma} H(\xi)=H(\gamma),
$$
for every $\gamma\in (0,\delta]$. Further, fix
\begin{eqnarray}\label{bound}
\lambda>\frac{\displaystyle\trace/2 +\sum_{i<j}a_{ij}}{\displaystyle\sum_{k=1}^{T}\alpha_k}\frac{\delta^2}{\displaystyle H(\delta)}.
\end{eqnarray}
By $(\rm h_2)$ there exists $\bar \gamma\in (0,\delta)$ such that
$$
h(t)<\frac{\lambda_1t}{\displaystyle\lambda\left(\sum_{k=1}^{T}\alpha_k\right)},
$$
for every $t\in ]0,\bar\gamma[$. Hence we obtain
$$
\sup_{t\in(0,\bar\gamma)}\frac{H(t)}{t^2}\leq \frac{\lambda_1}{2\displaystyle\lambda\left(\sum_{k=1}^{T}\alpha_k\right)}.
$$
On the other hand, bearing in mind \eqref{bound}, one has
$$
\left(\frac{\lambda_1}{\trace +2\sum_{i<j}a_{ij}}\right)\frac{\displaystyle\sum_{k=1}^{T}\alpha_k H(\delta)}{\delta^2}>\frac{\lambda_1}{2\lambda}.
$$

\noindent Taking
$$
0<\gamma<\min\left\{\bar \gamma,\displaystyle\left(\frac{\lambda_1}{\trace +2\sum_{i<j}a_{ij}}\right)^{-1/2}\delta\right\},
$$
and considering the above relations, one has
\begin{equation}\label{xy}
\frac{\displaystyle\left(\sum_{k=1}^{T}\alpha_k \right)\sup_{|\xi|\leq \gamma} H(\xi)}{\gamma^2}\leq \frac{\lambda_1}{2\lambda} <\left(\frac{\lambda_1}{\trace +2\sum_{i<j}a_{ij}}\right)\frac{\displaystyle\sum_{k=1}^{T}\alpha_k H(\delta)}{\delta^2}.
\end{equation}
\noindent Thus, condition $(\rm g_1)$ in Theorem \eqref{W2} is verified.\par

\noindent Now, we proceed by proving that condition ($\rm h_3$) implies ($\rm g_2$).
Indeed, by ($\rm h_3$), there are constants $\epsilon' \in \displaystyle\left]0,\lambda_1/\left(\sum_{k=1}^{T}\alpha_k\right)\right[,$ and
$\sigma>0$ such that
\begin{equation}\label{z}
{\frac{h(t)}{t}}<\displaystyle\frac{\lambda_1}{\displaystyle\sum_{k=1}^{T}\alpha_k}-\epsilon',
\end{equation}
 for every
$|t|\geq \sigma$. Since $h$ is a measurable locally bounded function, we also have
\begin{equation}\label{x}
M:=\sup_{t\in[-\sigma,\sigma]}|h(t)|<+\infty.
\end{equation}
\noindent Therefore, if $\xi\geq \sigma$, one has
$$
\int_0^\xi h(t)dt=\int_0^\sigma h(t)dt+\int_\sigma^\xi h(t)dt\leq M\sigma+\frac{1}{2}\left(\displaystyle {\lambda_1}/\left({\displaystyle\sum_{k=1}^{T}\alpha_k}\right)-\epsilon'\right)\xi^2,
$$
\noindent while, for $\xi\leq -\sigma$, it follows that
$$
\int_0^\xi h(t)dt=-\left[\int_\xi^{-\sigma} h(t)dt+\int_{-\sigma}^0 h(t)dt\right]\leq M\sigma+\frac{1}{2}\left(\displaystyle {\lambda_1}/\left({\displaystyle\sum_{k=1}^{T}\alpha_k}\right)-\epsilon'\right)\xi^2.
$$
\noindent Consequently,
\begin{eqnarray}\label{Dis}
\int_0^\xi h(t)dt\leq M\sigma+\frac{1}{2}\left(\displaystyle {\lambda_1}/\left({\displaystyle\sum_{k=1}^{T}\alpha_k}\right)-\epsilon'\right)\xi^2,\,\,\,\,\,\,\, \forall\, \xi\in\erre.
\end{eqnarray}
Hence, by using \eqref{Dis}, one has
$$\displaystyle\limsup_{|\xi|\rightarrow \infty}\frac{\displaystyle\int_0^\xi \alpha_k h(t)dt}{\xi^2}=\alpha_k\displaystyle\limsup_{|\xi|\rightarrow \infty}\frac{\displaystyle H(\xi)}{\xi^2}\leq\frac{1}{2}\displaystyle\left(\sum_{k=1}^{T}\alpha_k\right)\left(\displaystyle {\lambda_1}/\left({\displaystyle\sum_{k=1}^{T}\alpha_k}\right)-\epsilon'\right)<\frac{\lambda_1}{2}.$$
So, it is clear that condition ($\rm g_2$) holds, too. Finally, we achieve the
conclusion by applying Theorem \ref{Main}, taking into account that condition \eqref{xy} ensures that
$$
\lambda\in\left]\frac{\trace /2 +\sum_{i<j}a_{ij}}{\displaystyle \sum_{k=1}^{T}\alpha_k}\frac{\delta^2}{\displaystyle H(\delta)},\frac{\lambda_1}{2}\frac{\gamma^2}{\displaystyle\left(\sum_{k=1}^{T}\alpha_k \right)\sup_{|\xi|\leq \gamma} H(\xi)}\right[.
$$
Thus, our claim holds and the proof is completed.}
\end{pf}

\begin{remark}\rm{
As pointed out in the Introduction, in recent years, some related results on the existence of multiple solutions for algebraic systems of the form
\begin{eqnarray}\label{sistemino}
Au=\lambda g(u),
\end{eqnarray}
 were obtained by several authors (see the papers \cite{YZ,YZ2,WC}).\par
  \noindent We also emphasize for completeness that Theorem \ref{Main} and Corollary \ref{Main2} are mutually independent with respect to results and methods contained in the cited works also in the continuous setting, as can be easily seen by a direct computations.
}
\end{remark}

\begin{remark}\rm{
In the smooth case, Corollary 3.1 (and even the main abstract result) can be obtained directly from the celebrated Pucci-Serrin result (see \cite[Theorem 1.13]{KRV}); namely, for enough large parameters, a solution appears as a global minimum point for the energy functional (with negative energy level) while a second solution as a mountain pass one (with positive energy level). See, for instance, the recent paper \cite[Theorem 3.1]{MolicaRepovs} where, under suitable assumptions at zero and at infinity on the nonlinearity $g$, has been proved that, for every
$$\lambda>\left(\frac{\trace +2\sum_{i<j}a_{ij}}{2}\right)\left(\max_{\xi\neq
0}\frac{\sum_{k=1}^{T}G_k(\xi)}{\xi^2}\right)^{-1},$$
problem \eqref{sistemino} has at least
two distinct and nontrivial solutions.
}
\end{remark}

\section{Applications}
In this section we present some important cases of Theorem \ref{Main}. More precisely, we consider our previous result for some classes of discrete problems that appear in several technical applications.
\subsection{Tridiagonal inclusions.}

Let $T>1$ and $(a,b)\in
\erre^{-}\times\erre^{+}$ be such that
 \begin{equation}\label{first}
\cos\left(\frac{\pi}{T+1}\right)<-\frac{b}{2a}.
\end{equation}
\noindent Set
$$
{\T}_T(a,b,a)=\begin{pmatrix}
  b & a & 0 & ...  & 0 \\
  a & b & a & ...  & 0 \\
     &  & \ddots &  &  \\
  0 & ...  & a & b & a \\
  0 & ...  & 0 & a & b
\end{pmatrix}_{T\times T},
$$
\noindent and consider the following discrete problem
\begin{equation} \tag{$S_{\lambda}^{j}$} \label{tridiagonalcase}
L_{{\T}}(u)\in \lambda[j^-_k(u_k),j^+_k(u_k)],\quad\quad \forall\,k\in {\mathbb{Z}}[1,T],
\end{equation}
where
$$
L_{{\T}}(u):=\left\{
\begin{array}{l}
bu_1+au_2\\
au_{k-1}+bu_k+au_{k+1},\quad \forall\,k\in \{2,...,T-1\}\\
au_{T-1}+bu_T,
\end{array}
\right.
$$
the functions $j_k:\erre\rightarrow\erre$ are assumed to be measurable locally bounded, and $\lambda$ is a positive parameter. At this point, observing that ${\T}_T(a,b,a)$  is a symmetric and positive definite matrix whose first eigenvalue is given by
$$
\lambda_1=b+2a\cos\left(\frac{\pi}{T+1}\right),
$$
\noindent see
\cite[Example 9; p.179]{Sc}, we obtain the following multiplicity result.
\newpage
\begin{theorem}\label{tridiagonale}
Assume that
\begin{itemize}
\item[$(\rm g_1)$] There exist positive constants $\gamma$ and $\delta$, with $$\delta>\displaystyle\left(\frac{b+2a\cos\left({\pi}/{(T+1)}\right)}{bT+2a(T-1)}\right)^{1/2}\gamma,$$ such that
$$
\frac{\displaystyle\sum_{k=1}^{T}\sup_{|\xi|\leq \gamma} \displaystyle\int_0^\xi j_k(t)dt}{\gamma^2}<\frac{b+2a\cos\left({\pi}/{(T+1)}\right)}{bT+2a(T-1)}\frac{\displaystyle\sum_{k=1}^{T}\displaystyle\int_0^\delta j_k(t)dt}{\delta^2};
$$
\end{itemize}

\begin{itemize}
\item[$(\rm g_2)$] $\displaystyle\limsup_{|\xi|\rightarrow \infty}\frac{\displaystyle\int_0^\xi j_k(t)dt}{\xi^2}<\frac{b+2a\cos\left({\pi}/{(T+1)}\right)}{2},\,\,\,\,\,\,\forall\, k\in {\mathbb{Z}}[1,T].$

\end{itemize}
\noindent Then for every $\lambda$ belonging to
$$
\Lambda:=\left]\frac{bT+2a(T-1)}{2}\frac{\delta^2}{\displaystyle\sum_{k=1}^{T}\displaystyle\int_0^\delta j_k(t)dt},\frac{b+2a\cos\left({\pi}/{(T+1)}\right)}{2}\frac{\gamma^2}{\displaystyle \sum_{k=1}^{T}\sup_{|\xi|\leq \gamma} \displaystyle\int_0^\xi j_k(t)dt}\right[,
$$
problem \eqref{tridiagonalcase} admits at least three solutions.
\end{theorem}
The above result can be applied to second-order difference inclusions. Indeed,
it is well-know that the matrix
$${\T}_T(-1,2,-1)
:=\left(\begin{array}{ccccc}
  2 & -1 & 0 & ... & 0 \\
  -1 & 2 & -1 & ... & 0 \\
     &  & \ddots &  &  \\
  0 & ... & -1 & 2 & -1 \\
  0 & ... & 0 & -1 & 2
\end{array}\right)_{T\times T}\in {\mathfrak{X}}_{T},
$$
is associated to the second-order discrete boundary value problem
\begin{equation} \tag{$\widetilde{S}^j_\lambda$} \label{dxi}
\left\{
\begin{array}{l}
-\Delta^2u_{k-1}\in\lambda[j^-_k(u_k),j^+_k(u_k)],\quad \forall k  \in
{\mathbb{Z}}[1,T] \\
{u_0=u_{T+1}=0,}\\
\end{array}
\right.
\end{equation}
where $\Delta^2u_{k-1}:=\Delta(\Delta u_{k-1})$, and, as usual, $\Delta
u_{k-1}:=u_{k}-u_{k-1}$ denotes the forward difference operator.\par
\subsection{Fourth-order difference inclusions.} Boundary value problems involving fourth-order difference inclusions such as
\begin{equation} \tag{$D_{\lambda}^{g}$} \label{42order}
\left\{
\begin{array}{l}
-\Delta^4u_{k-2}\in \lambda[g^-_k(u_k),g^+_k(u_k)],\quad \forall\,k\in {\mathbb{Z}}[1,T]\\
{u_{-2}=u_{-1}=u_0=0}\\
u_{T+1}=u_{T+2}=u_{T+3}=0,\\
\end{array}
\right.
\end{equation}
can also be expressed as problem $(S_{A,\lambda}^{g})$, where $A$ is the real symmetric and positive definite matrix of the form
$$
A:=\left(
  \begin{array}{ccccccccc}
    6 & -4 & 1 & 0 & ... & 0 & 0 & 0 & 0 \\
    -4 & 6 & -4 & 1 & ... & 0 & 0 &0 & 0 \\
    1 & -4 & 6 & -4 & ... & 0 & 0 & 0 & 0 \\
    0 & 1 & -4 & 6 & ... & 0 & 0 & 0 & 0 \\
     &  &  &  & \ddots &  &  &  &  \\
    0 & 0 & 0 & 0 & ... & 6 & -4 & 1 & 0 \\
    0 & 0 & 0 & 0 & ... & -4 & 6 & -4 & 1 \\
    0 & 0 & 0 & 0 & ... & 1 & -4 & 6 & -4 \\
    0 & 0 & 0 & 0 & ... & 0 & 1 & -4 & 6 \\
  \end{array}
\right)\in {\mathfrak{X}}_{T}.
$$
In this case, let $\lambda_1$ be the first eigenvalue of $A$. Assuming that the following condition holds
\begin{itemize}
\item[$(\rm g_1')$] \textit{there exist two positive constants $\gamma$ and $\delta$, with $$\delta>\displaystyle\frac{\sqrt{\lambda_1}}{2}\gamma,$$ such that}
$$
\frac{\displaystyle\sum_{k=1}^{T}\sup_{|\xi|\leq \gamma} G_k(\xi)}{\gamma^2}<\frac{\lambda_1}{4}\frac{\displaystyle\sum_{k=1}^{T}G_k(\delta)}{\delta^2},
$$
\end{itemize}
in addition to $(\rm g_2)$, one has that, for every $\lambda$ belonging to
$$
\Lambda:=\left]\frac{2\delta^2}{\displaystyle\sum_{k=1}^{T} G_k(\delta)},\frac{\lambda_1}{2}\frac{\gamma^2}{\displaystyle\sum_{k=1}^{T}\sup_{|\xi|\leq \gamma} G_k(\xi)}\right[,
$$
problem \eqref{42order} admits at least three solutions.

\subsection{Partial difference inclusions.}
Nonlinear systems of the
form \eqref{W} arise in many applications as the boundary value problems for common boundary value problems
involving partial difference equations. For instance, here we consider the following problem, namely $(E_\lambda^f)$, given as follows
$$
4u(i,j)-u(i+1,j)-u(i-1,j)-u(i,j+1)-u(i,j-1)\in \lambda[f^-_{(i,j)}(u(i,j)),f^+_{(i,j)}(u(i,j))]
$$
for every $(i,j)\in {\mathbb{Z}}[1,m]\times{\mathbb{Z}}[1,n]$,
with boundary conditions
$$
u(i,0)=u(i,n+1)=0,\,\,\,\,\,\,\forall\, i\in {\mathbb{Z}}[1,m],
$$
$$
u(0,j)=u(m+1,j)=0,\,\,\,\,\,\,\forall\, j\in {\mathbb{Z}}[1,n],
$$
where every $f_{(i,j)}:\erre\rightarrow \erre$ denotes an essentially locally bounded function and $\lambda$ is a positive real parameter.\par
Let $z:{\mathbb{Z}}[1,m]\times{\mathbb{Z}}[1,n]\rightarrow {\mathbb{Z}}[1,mn]$ be the bijection defined by $z(i,j):=i+m(j-1),$ for every $(i,j)\in {\mathbb{Z}}[1,m]\times{\mathbb{Z}}[1,n]$.
 \noindent Let us denote $w_k:=u({z^{-1}(k)})$ and $g_k(w_k):=f_{z^{-1}(k)}(w_k)$, for every $k\in {\mathbb{Z}}[1,mn]$. With the above notations, problem $(E_\lambda^f)$ can be written as a nonlinear algebraic inclusion of the form
\begin{equation} \tag{$S_{B,\lambda}^{g}$} \label{W6}
\sum_{l=1}^{T}b_{kl}w_l\in \lambda[g^-_k(w_k),g^+_k(w_k)],\quad\quad \forall\,k\in {\mathbb{Z}}[1,mn],
\end{equation}
where $B$ is given by
$$
B:=(b_{ij})=
\left(
  \begin{array}{ccccccccc}
    L & -I_m & 0 & 0 & ... & 0 & 0 & 0 & 0 \\
    -I_m & L & -I_m & 0 & ... & 0 & 0 &0 & 0 \\
    0 & -I_m & L & -I_m & ... & 0 & 0 & 0 & 0 \\
    0 & 0 & -I_m & L & ... & 0 & 0 & 0 & 0 \\
     &  &  &  & \ddots &  &  &  &  \\
    0 & 0 & 0 & 0 & ... & L & -I_m & 0 & 0 \\
    0 & 0 & 0 & 0 & ... & -I_m & L & -I_m & 0 \\
    0 & 0 & 0 & 0 & ... & 0 & -I_m & L & -I_m \\
    0 & 0 & 0 & 0 & ... & 0 & 0 & -I_m & L \\
  \end{array}
\right)\in {\mathfrak{X}}_{mn},
$$
in which $L$ is the matrix of order $m$ defined by
$$
L:=
\left(
  \begin{array}{ccccccccc}
    4 & -1 & 0 & 0 & ... & 0 & 0 & 0 & 0 \\
    -1 & 4 & -1 & 0 & ... & 0 & 0 &0 & 0 \\
    0 & -1 & 4 & -1 & ... & 0 & 0 & 0 & 0 \\
    0 & 0 & -1 & 4 & ... & 0 & 0 & 0 & 0 \\
     &  &  &  & \ddots &  &  &  &  \\
    0 & 0 & 0 & 0 & ... & 4 & -1 & 0 & 0 \\
    0 & 0 & 0 & 0 & ... & -1 & 4 & -1 & 0 \\
    0 & 0 & 0 & 0 & ... & 0 & -1 & 4 & -1 \\
    0 & 0 & 0 & 0 & ... & 0 & 0 & -1 & 4 \\
  \end{array}
\right),
$$
\noindent and $I_m$ is the identity matrix. Denoting by $\lambda_{B}$ the first eigenvalue of $B$, we obtain the next result.
\begin{theorem}\label{Main3}
Assume that
\begin{itemize}
\item[$(\rm g_1)$] there exist two positive constants $\gamma$ and $\delta$, with $$\delta>\displaystyle\left(\frac{\lambda_B}{2(m+n)}\right)^{1/2}\gamma,$$ such that
$$
\frac{\displaystyle\sum_{k=1}^{mn}\sup_{|\xi|\leq \gamma} G_k(\xi)}{\gamma^2}<\frac{\lambda_B}{2(m+n)}\frac{\displaystyle\sum_{k=1}^{mn}G_k(\delta)}{\delta^2};
$$
\end{itemize}

\begin{itemize}
\item[$(\rm g_2)$] $\displaystyle\limsup_{|\xi|\rightarrow \infty}\frac{\displaystyle\int_0^\xi g_k(t)dt}{\xi^2}<\frac{\lambda_{B}}{2},\,\,\,\,\,\,\forall\, k\in {\mathbb{Z}}[1,mn].$

\end{itemize}
\noindent Then for every $\lambda$ belonging to
$$
\Lambda:=\left](m+n)\frac{\delta^2}{\displaystyle\sum_{k=1}^{mn} G_k(\delta)},\frac{\lambda_B}{2}\frac{\gamma^2}{\displaystyle\sum_{k=1}^{mn}\sup_{|\xi|\leq \gamma} G_k(\xi)}\right[,
$$
problem \eqref{W6} admits at least three solutions.
\end{theorem}
\begin{remark}\rm{
We observe that in \cite{JY}, Ji and Yang studied the structure of the spectrum of problem $(E_\lambda^f)$ by investigating the existence of a
positive eigenvector corresponding to the smallest positive eigenvalue $\lambda_B$. In conclusion, we refer to the paper of Galewski and Orpel \cite{GO} for some multiplicity results on discrete partial difference equations as well as to the monograph of Cheng \cite{Cheng} for their discrete geometrical interpretation.}
\end{remark}

\medskip
 \indent {\bf Acknowledgements.}  This paper was written when the first author was a visiting professor at the University of Ljubljana in 2012. He expresses his gratitude for the warm hospitality.
 The research was supported in part by the SRA grants P1-0292-0101 and J1-4144-0101.

\end{document}